\documentclass{article}
\usepackage{amssymb,amsmath,mathrsfs,latexsym,amscd}
\usepackage{exscale,latexsym,amsthm,graphics}

\usepackage{makeidx}
\makeindex
\usepackage{epsfig}
\usepackage{txfonts}
\newtheorem{thm}{Theorem}[section]
\newtheorem{cor}[thm]{Corollary}

\newtheorem{lemma}[thm]{Lemma}
\newtheorem{prop}[thm]{Proposition}

\newtheorem{remark}[thm]{Remark}

\def\bM {{\mathbb M}}

\def\qed{{\hfill $\Box$ \bigskip}}

\def\R {{\mathbb R}}

\def\P{{\mathbb P}}
\newcommand{\F}{\mathcal{F}}

\def\E{{\mathcal E}}
\def\P{{\mathbb P}}

\numberwithin{equation}{section}
\begin{document}

\noindent
{{\Large\bf Elliptic regularity results: n-regularized Liouville Brownian motion and non-symmetric diffusions associated with  degenerate forms}}

\bigskip
\noindent
{\bf Jiyong Shin}
\\

\noindent
{\small{\bf Abstract.} \\ 
We apply improved elliptic regularity results to a concrete symmetric Dirichlet form and various non-symmetric Dirichlet forms with possibly degenerate symmetric diffusion matrix. Given the (non)-symmetric Dirichlet form, using elliptic regularity results and stochastic calculus we show  weak existence of the corresponding singular stochastic differential equation for any starting point in some subset $E$ of $\R^d$. As an application of our approach we can show the existence of n-regularized Liouville Brownian motion only via Dirichlet form theory starting from all points in $\R^2$. \\



\section{Introduction}
In Dirichlet form theory, the general construction scheme of a Hunt process (associated to a given regular Dirichlet form) and the identification of the corresponding stochastic differential equation (SDE) may yield solutions starting from all points but exceptional set only. In recent years using elliptic regularity results based on \cite{BKR1}, \cite{BKR2}, and \cite[Theorem 1.7.4]{BKRS}, it has been shown for various (non)-symmetric Dirichlet forms  that one can obtain enough regularity of the corresponding semigroup of kernels and resolvent of kernels. In the sequel this leads to construct an associated diffusion process which solves the corresponding SDE (in the sense of martingale problem) starting from all explicitly specified points (see \cite{AKR}, \cite{BGS}). In \cite{AKR}, for instance, using the elliptic regularity results based on \cite{BKR1} and \cite{BKR2} the distorted Brownian motion associated to the given Dirichlet form is constructed, which solves martingale problem for explicitly specified set (i.e. $\{\rho >0\}$ there). Recently, by using the strict Fukushima decomposition the identification of the SDE  for the set $\{\rho >0\}$ in the sense of a weak solution has been worked out as a part of a general framework in \cite[Section 4]{ShTr13a}. 

In this paper we are concerned with applying well-known elliptic regularity results to more general (non)-symmetric Dirichlet forms, i.e. we consider the (non)-symmetric Dirichlet form (given as the closure of)
\begin{equation}\label{e;nsdde}
\E(f,g) : = \frac{1}{2} \int_{\R^d} \langle  A \nabla f , \nabla g \rangle \  \phi \ dx - \int_{\R^d} \langle B,\nabla f \rangle \ g \  \rho  dx, \quad
 f, g \in C_0^{\infty}(\R^d)
\end{equation}
on $L^2(\R^d,\rho dx)$. Here the conditions on $\rho$, $\phi>0$, (degenerate) symmetric diffusion matrix $A=(a_{ij})_{1 \le i,j \le d}$, and vector field $B : \R^d \to \R^d$ are formulated in (A) of Section \ref{symer} and (H1)-(H5) of Subsection \ref{s;ndee}, (H7) of Subsection \ref{s;3.2sdlce}, and (H8)-(H9) of Subsection \ref{s;3.3lms}.

The first aim of this paper is to take into account more general (non)-symmetric Dirichlet forms, i.e. $\phi \neq \rho$, the (degenerate) matrix $A$, and the vector field $B$ as in \eqref{e;nsdde}. In \cite{RoShTr} using elliptic regularity results in weighted spaces (\cite{BKR1} and \cite[Theorem 1.7.4]{BKRS}), stochastic calculus, and non-symmetric Dirichlet form theory, the (weak) existence of the non-symmetric distorted Brownian motion is shown for any starting point in some set in $\R^d$.
 Our extended result includes the consequences of \cite{RoShTr}, which is presented in Subection \ref{s;ndee} (see Remark \ref{r;ensdb}). In order to show the regularity of semigroup of kernels and resolvent kernels associated to the general Dirichlet forms as in \eqref{e;nsdde}, we adopt the improved elliptic regularity results as stated in \cite[Theorem 5.1] {BGS}) and Morrey's estimate \cite[Theorem 1.7.4]{BKRS}. This allows us to consider various reference measures, e.g. Lebesgue measure and different degeneracy of $A$ as in Section \ref{s;nonsydema}.

As a second aim of this paper we apply an elliptic regularity result to n-regularized Liouville Brownian motion. In \cite{GRV} C. Garban, R. Rhodes, V. Vargas constructed the Liouville Brownian motion (LBM). However, it is not clear as stated in \cite[Section 1.4]{GRV2} if one can construct LBM only via Dirichlet form theory.  We answer in the affirmative only for the case of n-regularized Liouville Brownian motion.  By setting $A=Id$, $B=0$, $\phi=1$, and $\rho > 0$ continuous on $\R^2$ we can obtain a weak solution to the corresponding SDE: 
\begin{equation}\label{e;elrlbn}
X_t = x+ \int_0^t  \rho^{-1/2} (X_s) \ dW_s, \quad t <  \zeta, \quad x \in \R^d,
\end{equation}
where $\zeta$ is the lifetime, $W = (W^1,\dots W^d)$ is a standard d-dimensional Brownian motion on $\R^d$. Applying this result to n-regularized Liouville Dirichlet form, we can conclude the existence of n-regularized Liouville Brownian motion only through Dirichlet form theory. 

The contents of this paper are organized as follows. In Section \ref{symer} we state some important elliptic regularity results corresponding to \eqref{e;nsdde} and consider one specific symmetric Dirichlet form.  We present elliptic regularity results and analytic consequences of the form. Subsequently, we then follow the methods and tools as in \cite{RoShTr} to show the weak existence of \eqref{e;elrlbn} for all $x \in \R^d$. As an application of this result, Subsection \ref{s;blbj} is devoted to the identification of  n-regualrized Liouville Brownian motion via Dirichlet form theory. In Section \ref{s;nonsydema} we consider 3 different types of (degenerate) non-symmetric Dirichlet forms of the form \eqref{e;nsdde}.

\section{A result of elliptic regularity and n-regularized Liouville Brownina motion}\label{symer}
{\bf{Notations}}:\\
For open set $U \subset \R^d$, $d \ge 2$ with Borel  $\sigma$-algebra $\mathcal{B}(U)$, we denote the set of all $\mathcal{B}(U)$-measurable $f : U \rightarrow \R$ which are bounded, or nonnegative by $\mathcal{B}_b(U)$, $\mathcal{B}^{+}(U)$, respectively. The usual $L^q$-spaces $L^q(U, \mu)$, $q \in[1,\infty]$ are equipped with $L^{q}$-norm $\| \cdot \|_{L^q(U,\mu)}$ with respect to the  measure $\mu$ on $U$ and $\mathcal{A}_b$ : = $\mathcal{A} \cap \mathcal{B}_b(U)$ for $\mathcal{A} \subset L^q(U,\mu)$,  and $L^{q}_{loc}(U,\mu) := \{ f \,|\; f \cdot 1_{U} \in L^q(U, \mu),\,\forall  G \subset U, G \text{ relatively compact open} \}$, where $1_A$ denotes the indicator function of a set $A$. The inner product on $L^2(E, \mu)$ is denoted by $(\cdot,\cdot)_{L^2(E, \mu)}$. Let $\nabla f : = ( \partial_{1} f, \dots , \partial_{d} f )$  and  $\Delta f : = \sum_{j=1}^{d} \partial_{jj} f$ where $\partial_j f$ is the $j$-th weak partial derivative of $f$ and $\partial_{ij} f := \partial_{i}(\partial_{j} f) $, $i,j=1, \dots, d$.  As usual $dx$ denotes Lebesgue measure on $\R^d$ and the Sobolev space $H^{1,q}(U, dx)$, $q \ge 1$ is defined to be the set of all functions $f \in L^{q}(U, dx)$ such that $\partial_{j} f \in L^{q}(U, dx)$, $j=1, \dots, d$, and 
$H^{1,q}_{loc}(U, dx) : =  \{ f  \,|\;  f \cdot \varphi \in H^{1,q}(U, dx),\,\forall \varphi \in  C_0^{\infty}(U)\}$. 
Here $C_0^{\infty}(U)$ denotes the set of all infinitely differentiable functions with compact support in $U$. We also denote the set of continuous functions on $U$, the set of continuous bounded functions on $U$, the set of compactly supported continuous functions in $U$ by $C(U)$, $C_b(U)$, $C_0(U)$, respectively. $C_{\infty}(U)$ denotes the space of continuous functions on $U$ which vanish at infinity and  $C^{1-\alpha}_{loc}(U)$, $0<\alpha<1$, denotes the set of all locally H\"{o}lder continuous functions of order $1-\alpha$ on $U$. We equip $\R^d$ with the Euclidean norm $\| \cdot \|$, the corresponding inner product $\langle \cdot, \cdot \rangle$ and write $B_{r}(x): = \{ y \in \R^d \ | \ \|x-y\| < r  \}$, $x \in \R^d$.\\

In this section we present a result of symmetric elliptic regularity and an application to n-regularized Liouville Brownian motion. 
\subsection{An elliptic regularity result of symmetric Dirichlet forms}
 Throughout this subsection, we assume
\begin{itemize}
\item[(A)]
$\rho(x)>0$, $\forall x \in \R^d$ and $\rho$ is continuous on $\R^d$.
\end{itemize}
We consider the symmetric positive definite bilinear form
\[
\E(f,g) : = \frac{1}{2} \int_{\R^d} \langle \nabla f ,  \nabla g \rangle \ dx , \quad f, g \in C_0^{\infty}(\R^d).
\]
Then $(\E,C_0^{\infty}(\R^d))$ is closable in $L^2(\R^d,m)$, $m:=\rho dx$ and its closure $(\E,D(\E))$ is a symmetric, strongly local, regular Dirichlet form (cf. \cite[II. Section 2. a)]{MR}). Furthermore,  it is well known that $(\E,C_0^{\infty}(\R^d))$ is closable in $L^2(\R^d,dx)$ and its closure is the symmetric Dirichlet form $(\E, H^{1,2}(\R^d,dx))$. Let us state the properties of $(\E,D(\E))$ by using the well-known Dirichlet form  $(\E, H^{1,2}(\R^d,dx))$:

\begin{thm}\label{t;retrcos}
\begin{itemize}
\item[(i)]  The Dirichlet form $(\E,D(\E))$ is the trace Dirichlet form of $(\E,H^{1,2}(\R^d,dx))$ relative to $m$ and $D(\E) = \{ f \in L^2(\R^d, m) \cap H^{1,2}_{loc} (\R^d, dx) \ | \ \nabla f \in L^2(\R^d, dx) \}$. 
\item[(ii)] Let $d=2$. Then the Dirichlet form $(\E,D(\E))$ is recurrent. In particular, $(\E,D(\E))$ is conservative.
\item[(iii)] Let $d \ge 3$.  Then the Dirichlet form $(\E,D(\E))$ is transient.
\end{itemize}
\end{thm} 
\proof
Recall that there exists a Brownian motion $((W_t)_{t \ge 0}, (\P_x)_{x \in \R^d})$ on $\R^d$ associated with the Dirichlet form $(\E,H^{1,2}(\R^d,dx))$ (cf. \cite[Example 4.2.1]{FOT}). The positive Radon measure $m =\rho dx$ does not charge capacity zero set of $(\E, H^{1,2}(\R^d,dx))$ and the positive continuous additive functional $(A_t)_{t \ge 0}$ (in the strict sense) associated with the Revuz measure $m$ is given by (cf. \cite[Example 5.1.1]{FOT})
\[
A_t = \int_0^t \rho(W_s) \ ds, \quad t \ge 0.
\]
 We define the support of  $(A_t)_{t \ge 0}$ by
\[
\tilde{Y} = \{ x \in \R^d \ | \ \P_x (R=0) =1\},
\]
where $R = \inf\{ t > 0 \ | \ A_t > 0 \}$. Clearly, the support of $m$ and $(A_t)_{t \ge 0}$ is $\R^d$.
Therefore by \cite[Section 6.2]{FOT}, the trace Dirichlet form of $(\E, H^{1,2}(\R^d,dx))$ on $\R^d$ relatve to $m$ is $(\E,\mathcal{F})$, where
\[
\mathcal{F} = \{ f \in L^2(\R^d, m)  \ | \  f \in H_e^{1,2}(\R^d,dx)  \}.
\]
By \cite[Theorem 2.2.13]{ChFu}, $H_e^{1,2}(\R^d,dx) = \{ f \in L^2_{loc}(\R^d, dx)  \ | \  \nabla f \in  L^2(\R^d,dx)  \}$. Here we denote $H_e^{1,2}(\R^d,dx)$ by the extended Dirichlet space of $H^{1,2}(\R^d,dx)$ (see \cite[p. 40]{FOT} for the definition). Therefore
\[
\mathcal{F} = \{ f \in L^2(\R^d, m) \cap H^{1,2}_{loc} (\R^d, dx) \ | \ \nabla f \in L^2(\R^d, dx) \}.
\]
Furthermore, by \cite[Theorem 6.2.1]{FOT} the core of $(\E,\mathcal{F})$ is $C_0^{\infty}(\R^d)$. This implies that the trace Dirichlet form of $(\E, H^{1,2}(\R^d,dx))$ relative to the measure $m$ is the Dirichlet form $(\E,D(\E))$ (see \cite[Section 6.2]{FOT}). Note that the 2-dimensional (resp. d-dimensional, $d \ge 3$) Brownian motion $(W_t)_{t \ge 0}$ is recurrent on $\R^2$ (resp. transient on $\R^d$, $d \ge 3$) and therefore by \cite[Theorem 6.2.3 (ii)]{FOT}, $(\E,D(\E))$ is recurrent (resp. transient).
\qed

Let  $(T_t)_{t > 0}$ and $(G_{\alpha})_{\alpha > 0}$ be the $L^2(\R^d, m)$-semigroup  and resolvent associated to $(\E,D(\E))$ and $(L,D(L))$ be the corresponding generator. $(T_t)_{t >0}$ and $(G_{\alpha})_{\alpha > 0}$ restricted to $L^1(\R^d,m) \cap L^{\infty}(\R^d,m)$ can be extended to strongly continuous contraction semigroups  on all $L^r(\R^d,m)$, $r \in [1,\infty)$ and we denote the corresponding operator families again by $(T_t)_{t > 0}$ and $(G_{\alpha})_{\alpha > 0}$. Let $(L_r, D(L_r))$ be the corresponding generator on $L^r(\R^d,m)$, $r \in [1,\infty)$. Furthermore  $(T_t)_{t > 0}$ is an analytic semigroup on every $L^r(\R^d,m)$, $r \in (1,\infty)$ (cf. \cite[Proposition 1.8, Remark 1.2]{LiSe}).    It is easy to see that 
$C_0^{\infty} (\R^d) \subset D(L_r)$ and for all $\varphi \in C_0^{\infty} (\R^d)$
\begin{equation}\label{eq;slbmg}
L \varphi = \frac{1}{2} \ \Delta \varphi \cdot \rho^{-1}.
\end{equation}

Let us first restate an improved elliptic regularity result as stated in \cite[Theorem 5.1] {BGS}  to prove regularity of $(T_t)_{t>0}$ and $(G_{\alpha})_{\alpha > 0}$ (cf. \cite{BKR2}):
\begin{prop}\label{p;gdegep}
Let $U \subset \R^d$, $d \ge 2$ be an open and $\mu$ a locally finite (signed) Borel measure on $U$ that is absolutely continuous w.r.t. $dx$ on $U$. Suppose $d_{ij} \in H^{1,p}_{loc}(U,dx)$ for some $p > d$ and the matrix $D= (d_{ij})_{1 \le i,j \le d}$ is locally strictly elliptic $dx$-a.e. on $U$, i.e. for each compact set $K \subset U$, there exists some $\kappa_K > 0$ such that  $\kappa_K \| \xi \|^2 \le  \langle D(x) \xi, \xi \rangle$, $\forall \xi \in \R^d$ $dx$-a.e. on $K$. Let either $h_i$, $c \in L_{loc}^p(U,dx)$ or $h_i$, $c \in L^p_{loc}(U, \mu)$ and let $f \in L^p_{loc}(U,dx)$. Assume that one has
\[
\int_{U} \Big( \sum_{i,j = 1}^{d} d_{ij} \partial_{ij} \varphi + \sum_{i=1}^{d} h_i \partial_i \varphi + c \varphi  \Big) \  d\mu = \int_U  \varphi \ f \ dx, \quad \forall \varphi \in C_0^{\infty} (U),
\]
where we assume that $h_i$, $c$ are locally $\mu$-integrable. Then $\mu$ has a density in $H^{1,p}_{loc}(U)$ that is locally H\"older continuous.
\end{prop}
Additionally, we restate the Morrey's estimate in our setting  (see \cite[Theorem 1.7.4]{BKRS}).
\begin{prop}\label{t;morrey2.4}
Assume $p > d \ge 2$. Let V be a bounded domain in $\R^d$ and $b : V \to \R^d$ and $c,e : V \to \R$ such that
\[
h_i \in L^p (V, dx) \ \  \text{and} \ \ c,e \in L^q(V,dx) \ \ \text{for} \ \ q:=\frac{dp}{d+p}>1.
\]
Let $a_{ij} = a_{ji}$, $a_{ij} \in H^{1,p}_{loc}(\R^d,dx) $ for all $1 \le i,j \le d$ and $\kappa^{-1} \ \| \xi \|^2 \le \langle A(x) \xi,\xi \rangle \le \kappa \| \xi \|^2$, $\forall \xi \in \R^d$, $x \in V$ for some $\kappa \ge 1$. Assume that $u \in H^{1,p}(V)$ is a solution of 
\[
\int_{V} \sum_{i=1}^d \Big( \partial_i \varphi \Big( \sum_{j=1}^d a_{ij} \partial_j u + h_i u  \Big)  \Big) + \varphi (cu + e) \ dx =0, \quad \forall \varphi \in C_0^{\infty}(V),
\]
Then for every domain $V'$ with $V' \subset \overline{V'} \subset V$, we obtain the estimate
\[
\| u \|_{H^{1,p}(V')} \le c (\| e \|_{L^q(V,dx)} + \| u \|_{L^1(V,dx)}),
\]
where  $c < \infty$ is some constant independent of $e$ and $u$.
\end{prop}

\begin{cor}\label{co;resyg}
Let $\alpha > 0$.
Suppose $g \in L^r(\R^d, m)$, $r \in [p,\infty)$. Then
\[
G_{\alpha} g \in H_{loc}^{1,p}(\R^d,dx)
\]
and for any open ball $B'\subset \overline{B'}\subset B\subset \R^d$   there exists $c_{B,\alpha} \in (0,\infty)$, independent of $g$, such that 
\begin{equation}\label{sonosyr}
\|  G_{\alpha} g \ \|_{H^{1,p}(B',dx)} \le c_{B,\alpha} \ \Big(\|G_{\alpha} g \|_{L^1(B,m)} + \| g \|_{L^{p}(B,m)} \Big).
\end{equation}
\end{cor}
\proof
Let $g \in C_0^{\infty}(\R^d)$. Then we have
\begin{equation}\label{selrqe}
\int_{\R^d} (\alpha - L  ) \varphi \ G_{\alpha} g \ \rho \ dx = \int_{\R^d} \varphi \ g \ \rho \ dx, \quad \forall \varphi \in C_0^{\infty}(\R^d).
\end{equation}
Since $G_{\alpha} g \in L^1(\R^d,m) \cap L^{\infty} (\R^d,m)$, we can define the signed Radon measure $\mu = - \frac{1}{2} G_{\alpha}g dx$. Now we apply Proposition \ref{p;gdegep} with $\mu = - \frac{1}{2} G_{\alpha} g dx$, $D=Id$, $c : = - 2 \alpha \ \rho \ \varphi$ and $f=g \rho$ to prove the assertion for $g \in C_0^{\infty}(\R^d)$. Clearly, all integrability conditions are satisfied. Therefore, $G_{\alpha} g \in H_{loc}^{1,p}(\R^d,dx)$ for $g \in C_0^{\infty}(\R^d)$. Let $u:= \frac{1}{2} G_{\alpha} g$. Using integration by parts \eqref{selrqe} can be rewritten as
\[
\int_{\R^d} \sum_{i=1}^{d} \partial_i \varphi \ \partial_i u + \varphi (2 \alpha \rho u - g \rho) \ dx = 0, \quad \forall \varphi \in C_0^{\infty}(\R^d).
\]
Then for any open balls $B'$, $B$ with $B'\subset \overline{B'}\subset B \subset \R^d$  we apply Proposition \ref{t;morrey2.4} with $\varphi \in C_0^{\infty}(B)$, $c = 2 \alpha \rho \in L^q(B,dx)$ and $e= -\rho g \in L^q(B,dx)$, $q:= dp/(d+p) > 1$. Hence we obtain \eqref{sonosyr} for $g \in C_0^{\infty}(\R^d)$. Since $C_0^{\infty}(\R^d)$ is dense in ($L^r(\R^d, m)$, $\| \cdot \|_{L^r(\R^d, m)})$, $r \in [p,\infty) $, the assertion for general $g \in L^r(\R^d, m)$ follows by continuity.
\qed

\begin{cor}\label{c;stranw}
Let $t>0$, $r \in [p,\infty)$. 
\begin{itemize}
\item[(i)] Let $u \in D(L_r)$. Then 
\[
T_t u \in H^{1,p}_{loc} (\R^d, dx)
\]
and for any open ball $B'\subset \overline{B'}\subset B \subset \R^d$ there exists $c_B \in (0,\infty)$ (independent of $u$ and $t$) such that
\begin{equation}\label{eq;sylbm}
\| T_t u\|_{H^{1,p}(B',dx)} \le c_B \left( m(B)^{\frac{r-1}{r}} \| u \|_{L^r (\R^d,m) }  + m(B)^{\frac{r-p}{rp}} \| (1-L_r) u \|_{L^r(\R^d,m)}    \right).
\end{equation}
\item[(ii)] Let $u \in L^r(\R^d,m)$. Then the above statements still hold with \eqref{eq;sylbm} replaced by
\[
\|T_t u\|_{H^{1,p}(B',dx)} \le  \tilde{c}_B \left(1+ \frac{1}{t}\right)  \| u \|_{L^r(\R^d,m)},
\]
where $\tilde{c}_B \in (0,\infty)$ (independent of $f$, $t$).
\end{itemize}
\end{cor}
\proof
Using Corollary \ref{co;resyg} and the analyticity of $(T_t)_{t > 0}$ one can show (i) and (ii) (see \cite[Corollary 2.4]{AKR}).
\qed
\begin{remark}
The assumption that continuous $\rho$ is strict positive on $\R^d$ as in (A) plays a crucial role in the proof of  \eqref{sonosyr} and \eqref{eq;sylbm}. Therefore we can not relax the strict positivity of $\rho$.
\end{remark}

By virtue of Corollary \ref{co;resyg} and Corollary \ref{c;stranw}, exactly as in \cite[Section 3]{AKR}, we obtain a semigroup of kernels $(P_t)_{t > 0}$ and resolvent of kernels $(R_{\alpha})_{\alpha >0}$ on $\R^d$ as follows:\\
There exists a transition kernel density $p_t(\cdot,\cdot)$ on $\R^d$
such that
\[
P_t f(x) : = \int_{\R^d} f(y) p_t(x,y) \ m(dy), \quad x \in \R^d, \ t > 0
\]
is a (temporally homogeneous) sub-Markovian transition function (cf. \cite[1.2]{CW}) and an $m$-version of $T_t f $ for any $f \in \cup_{r \ge p} L^r(\R^d,m)$. Moreover, letting $P_0 : = id$ and following \cite[Proposition 3.2]{AKR}  we obtain:

\begin{thm}\label{t;2.7gp}
\begin{itemize}
\item[(i)] $(P_t)_{t>0}$ is a semigroup of kernels on $\R^d$ which is $L^r(\R^d,m)$-strong Feller for all $r \in [p,\infty)$, i.e.
\[
P_t f \in C(\R^d), \quad \forall f \in \cup_{r \ge p} L^r(\R^d,m), \quad \forall t > 0.
\]
\item[(ii)]
\[
\lim_{t \to 0} P_{t+s} f(x) = P_s f(x), \quad \forall s \ge 0, \ x \in \R^d, \ f \in C_0^{\infty}(\R^d).
\]
\item[(iii)]
$(P_t)_{t > 0}$ is a measurable semigroup on $\R^d$, i.e. for $f \in \mathcal{B}^+ (\R^d)$ the map $(t,x) \mapsto P_t f(x)$ is $\mathcal{B}([0,\infty) \times \R^d)$-measurable.

\end{itemize}
\end{thm}

Similarly, since for $\alpha > 0$, $f \in \cup_{r \ge p} L^r(\R^d,m)$, $G_{\alpha} f$ has a unique continuous $m$-version on $\R^d$ by Corollary \ref{co;resyg} as in \cite[Lemma 3.4, Proposition 3.5]{AKR}, we can find $(R_{\alpha})_{\alpha > 0}$ with resolvent kernel density $r_{\alpha}(\cdot, \cdot)$ defined on $\R^d \times \R^d$ such that
\[
R_{\alpha} f(x) : = \int f(y)\, r_{\alpha} (x,y) \ m(dy),  \quad x \in \R^d, \ \alpha > 0,
\]
satisfies 
\[
R_{\alpha}f  \in C(\R^d) \ \text{and} \ R_{\alpha} f = G_{\alpha} f  \ \  m\text{-a.e for any} \  f \in \cup_{r \ge p} L^r(\R^d,m).
\]
Suppose that $(\E,D(\E))$ is conservative  (e.g. Theorem \ref{t;retrcos} (ii)). Then we obtain from \cite[Proposition 3.8]{AKR} the following:
\begin{thm}\label{t;trsyst}
Suppose that $(\E,D(\E))$ is conservative (e.g. d=2). Then,
\begin{itemize}
\item[(i)] $\alpha R_{\alpha} 1(x) = 1$ for all $x \in \R^d$, $\alpha > 0$.
\item[(ii)] $(P_t)_{t>0}$ is strong Feller on $\R^d$, i.e. $P_t(\mathcal{B}_b(\R^d)) \subset C_b(\R^d)$ for all $t>0$.
\item[(iii)] $P_t 1(x) = 1$ for all $x \in \R^d$, $t>0$.
\end{itemize}
\end{thm}

\begin{thm}\label{th;ehplb}
There exists a continuous Hunt process (i.e. strong Markov with continuous sample path)
\[
\bM =  (\Omega, \F, (\F_t)_{t \ge 0}, \zeta, (X_t)_{t \ge 0}, (\P_x)_{x \in \R^d \cup \{\Delta\}  }  )
\]
with state space $\R^d$, cemetery $\Delta$, and the lifetime $\zeta:=\inf\{ t \ge 0 \ | \ X_t=\Delta\}$, having the transition function $(P_t)_{t \ge 0}$ as transition semigroup.
\end{thm}
\proof
Following the proof of \cite[Section 4]{AKR} with the state space $\{ \rho > 0 \}$ replaced by $\R^d$, we can construct the associated diffusion process having the transition function $(P_t)_{t \ge 0}$ as transition semigroup on the state space $\R^d$.
\qed

\begin{remark}\label{r;sytr}
Suppose that $(\E,D(\E))$ is conservative (e.g. d=2). Then, by Theorem \ref{t;trsyst} (ii), $\bM$ becomes a classical (conservative) diffusion with state space $\R^d$, i.e.
\[
\P_x(\zeta = \infty) = 1, \quad \forall x \in \R^d.
\]
\end{remark}

\begin{lemma}\label{l;trsydepsr}
\begin{itemize}
\item[(i)] Let $f \in \bigcup_{r \in [p, \infty)} L^r(\R^d,m)$, $f \ge 0$, then for all $t > 0$, $x \in \R^d$,
\[
\int_0^t P_s f(x) \ ds < \infty,
\]
hence
\[
\int \int_0^t f(X_s) \ ds \ d\P_x < \infty.
\]
\item[(ii)] Let $u \in C_0^{\infty} (\R^d)$, $\alpha > 0$. Then
\[
R_{\alpha} \big( (\alpha - L) u\big)(x) = u(x) \quad \forall x \in \R^d. 
\]
\item[(iii)] Let $u \in C_0^{\infty}(\R^d)$, $t>0$. Then
\[
P_t u(x) - u(x) = \int_0^t P_s(Lu)(x) \ ds \quad \forall x \in \R^d.
\]
\end{itemize}
\end{lemma}
\proof
The proof is the same as \cite[Lemma 5.1]{AKR} with $E$ replaced by $\R^d$. 
\qed

\begin{lemma}\label{l;sytrma}
\begin{itemize}
\item[(i)]
For $u \in C_0^{\infty} (\R^d)$
\[
L u^2 - 2 u \ Lu =   \| \nabla u \|^2 \cdot \rho^{-1}.
\]
\item[(ii)]
Let $u \in C_0^{\infty}(\R^d)$ and
\[
M_t : = \left( u(X_t) - u(X_0) - \int_0^t Lu(X_r) \ dr \right)^2 - \int_0^t  \Big(   \| \nabla u \|^2 \cdot \rho^{-1} \Big) (X_r) \ dr, \quad t \ge 0.
\]
Then $(M_t)_{t \ge 0}$ is  an $(\mathcal{F}_t)_{t \ge 0}$-martingale under $\P_x$, $\forall x \in \R^d$.

\end{itemize}
\end{lemma}
\proof
(i) By \eqref{eq;slbmg} we obtain for $u \in C_0^{\infty} (\R^d)$
\[
L u^2 - 2 u \ Lu = \| \nabla u \|^2 \ \rho^{-1} + u \Delta u \ \rho^{-1} - u  \Delta u \ \rho^{-1}  = \|\nabla u\|^2 \cdot \rho^{-1}.
\]
(ii) Using Lemma \ref{l;trsydepsr} and (i), the proof is similar to \cite[Proposition 3.3]{RoShTr} with $E$ replaced by $\R^d$ (cf.  \cite[Section 6]{RoShTr}  ).
\qed

\begin{lemma}\label{potwitrstr}
Let $(B_k)_{k \ge 1}$ be an increasing sequence of relatively compact open sets in $\R^d$ with $\cup_{k \ge 1} B_k= \R^d$.
Then
for all $x \in \R^d$
\[
\P_x \Big(\lim_{k \rightarrow \infty} \sigma_{\R^d \setminus B_{k}} \ge \zeta \Big)=1.
\] 
\end{lemma}
\proof
The proof is the same as \cite[Lemma 3.4]{RoShTr} with $E$ replaced by $\R^d$. 
\qed

\begin{thm}\label{t;lbmdf}
For each $x\in \R^d$, the process $\bM$ satisfies 
\[
X_t = x+ \int_0^t  \rho^{-1/2} (X_s) \ dW_s, \quad t < \zeta,
\]
$\P_x$-a.s. where $W = (W^1,\dots,W^d)$ is a standard d-dimensional Brownian motion on $\R^d$. In particular, if $(\E,D(\E))$ is conservative (e.g. $d=2$), then the lifetime $\zeta$ is replaced by $\infty$.
\end{thm}
\proof
Using Lemma \ref{l;sytrma} and Lemma \ref{potwitrstr}, the proof is similar to the proof of Theorem \ref{t;fdnons}. So we omit it.
\qed

\begin{remark}
In addition  to (A), suppose that $\rho \in H^{1,1}_{loc} (\R^d,dx)$ and $\| \nabla (  \sqrt{\rho}^{-1} )  \| \in L^{2(d+1)}_{loc} (\R^d,dx)$. Then by \cite[Theorem 1.1]{Zh} the solution in Theorem \ref{t;lbmdf} is strong, pathwise unique, and non explosive.

\end{remark}

\subsection{Application to n-regularized Liouville Brownian motion}\label{s;blbj}
In \cite{GRV} C. Garban, R. Rhodes, V. Vargas constructed a diffusion process with a massive Gaussian free field X, called Liouville Brownian motion (LBM). By classical theory of Gaussian multiplicative chaos (cf. \cite{Kah}) Liouville measure is well defined and subsequently the LBM is constructed  by using approximation and construction of positive continuous additive functional of Liouville measure relevant to a Brownian motion on $\R^2$. However as stated in \cite[Section 1.4]{GRV2}  it is not clear if one can construct directly the Liouville Brownian motion via the theory of Dirichlet forms without using the results in \cite{GRV}. In general there is no theory of  Dirichlet forms which enables to get rid of the polar set and construct a Hunt process starting from all points of $\R^2$. We can answer in the affirmative only for the case of n-regularized Liouville Brownian motion (n-LBM). Let us first restate some definitions and results as stated in \cite{GRV} and \cite{GRV2}:\\

\noindent  {\bf{Massive Gaussian free field:}}\\
Given a real number $m > 0$, we consider a whole plane massive Gaussian free field which is a centered Gaussian random distribution (in the sense of Schwartz)  on a probability space $(\Omega, \mathcal{A}, P)$ with covariance function given by the Green function $G^{(m)}$ of the operator $m^2 -\Delta$, i.e.
\[
(m^2 - \Delta) G^{(m)}(x, \cdot) = 2 \pi \delta_x, \quad x \in \R^2,
\] 
where $\delta_x$ stands for the Dirac mass at $x$.
The massive Green function with the operator $(m^2 - \Delta)$ can be written as 
\[
G^{(m)}(x,y) = \int_0^{\infty} e^{- \frac{m^2}{2}s - \frac{\| x-y \|^2}{2s}} \frac{ds}{2s} 
= \int_{1}^{\infty} \frac{k_{m} (s (x-y)) }{s} \ ds, \quad x,y \in \R^2,
\]
where 
\[
k_m (z) = \frac{1}{2} \int_0^{\infty} e^{- \frac{m^2}{2s} \| z  \|^2 - \frac{s}{2} } \ ds.
\]
Note that this massive Green function is  a kernel of $\sigma$-positive type in the sense of Kahane \cite{Kah} since we integrate a continuous function of positive type w.r.t. a positive measure.
Let $(c_n)_{n \ge 1}$ be an unbounded strictly increasing sequence such that $c_1 = 1$ and  $(Y_n)_{n \ge 1}$ be a family of  independent centered continuous Gaussian fields on $\R^2$ on the probability space $(\Omega, \mathcal{A}, P)$ with covariance kernel given by 
\[
E[Y_n (x) \ Y_n (y)] = \int_{c_{n-1}}^{c_n}  \frac{k_m (s(x-y))}{s} ds.
\]
The massive Gaussian free field is the Gaussian distribution defined by 
\[
X(x) = \sum_{k \ge 1} Y_k (x).
\]
\noindent  {\bf{Liouville Brownian motion:}}\\
We define $n$-regularized field by
\[
X_n(x) = \sum_{k=1}^{n} Y_k(x), \quad n \ge 1
\]
and the associated  $n$-regularized Liouville measure by
\[
M_{n,\gamma} (dz)=  \exp \Big( \gamma X_n(z) - \frac{\gamma^2}{2} E[X_n(z)^2] \Big) \ dz, \quad \gamma \in (0,2).
\]
Since $M_{n,\gamma}$ does not charge any polar set of the given planar Brownian motion $(W_t)_{t \ge 0}$ (independent of $X$), $P$-a.s. we can deduce a unique positive continuous additive functional $(F_t^n)_{t \ge 0}$ in the strict sense associated to $M_{n,\gamma}$ (see \cite{FOT}). For every $n \ge 1 $, $F_t^n$ is defined as
\[
F_t^n: =\int_0^t  e^{\gamma  X_n(W_s) - \frac{\gamma^2}{2} E[X_n(W_s)^2]   }  ds,  \quad t \ge 0, 
\]
which is strictly increasing in $t$ and then the n-LBM is defined as
\[
\mathcal{B}^n_t = W_{(F_t^{n})^{-1}}, \quad t \ge 0.
\]
The n-LBM $(\mathcal{B}^n_t)_{t \ge 0}$ can be thought of as  the solution of the stochastic differential equation:
\[
d \mathcal{B}^n_t = e^{-\frac{\gamma}{2} X_n(\mathcal{B}^n_t) + \frac{\gamma^2}{4} E[X_n(\mathcal{B}^n_t)^2]   } dB_t, \quad \mathcal{B}^n_0 =x, \quad x \in \R^2,
\]
where $B$ is a standard Brownian motion on $\R^2$ independent of $X$. Furthermore n-LBM is associated to the n-regularized Dirichlet form
\[
\E^n(f,g) = \frac{1}{2}  \int_{\R^2} \nabla f (x) \cdot \nabla g(x) \ dx, \quad f,g \in D(\E^n),
\]
where $D(\E^n) = \{ f \in L^2(\R^2, M_{n,\gamma}) \cap H^{1,2}_{loc} (\R^2, dx) \ | \ \nabla f \in L^2(\R^d, dx) \}$.\\
By the classical theory of Gaussian multiplicative chaos (see \cite{Kah}), $P$-a.s. the family $(M_{n,\gamma})_{n \ge 1}$ weakly converges to a Radon measure $M_{\gamma}$, which is called Liouville measure, i.e. the Liouville measure $M_{\gamma}$ is formally defined as 
\[
M_{\gamma} (dz)=  \exp \Big( \gamma X(z) - \frac{\gamma^2}{2} E[X(z)^2] \Big) \ dz.
\]
In \cite{GRV} the positive continuous additive functional $(F_t)_{t \ge 0}$ (in the strict sense) of the planar Brownian motion $(W_t)_{t \ge 0}$ associated with the measure $M_{\gamma}$ is constructed  and then the LBM is defined as
\[
\mathcal{B}_t = W_{F_t^{-1}}, \quad t \ge 0.
\]
The support of the positive continuous additive functional $(F_t)_{t \ge 0}$ is the whole space $\R^2$  (see \cite{GRV}). Since the support of $(F_t)_{t \ge 0}$ is the whole space $\R^2$, the LBM is associated to the Liouville Dirichlet form (see \cite[Theorem 6.2.1]{FOT})
\[
\E(f,g) = \frac{1}{2}  \int_{\R^2} \nabla f (x) \cdot \nabla g(x) \ dx, \quad f,g \in D(\E),
\]
where $D(\E) = \{ f \in L^2(\R^2, M_{\gamma}) \cap H^{1,2}_{loc} (\R^2, dx) \ | \ \nabla f \in L^2(\R^d, dx) \}$.\\

\noindent  {\bf{Construction of n-regularized Liouville Brownian motion via Dirichlet form theory:}}\\
 Using the elliptic regularity results stated in this section, we construct the n-LBM starting from all points on $\R^2$ via Dirichlet form theory. Let
\[
 \rho (z): = e^{ \gamma X_n(z) - \frac{\gamma^2}{2} E[X^2_n(z)] }, \quad z \in \R^2. 
\]
Then, clealry $\rho$ satisfies the assumption (A). Now we consider the symmetric bilinear form
\[
\E^n (f,g) : = \frac{1}{2}\int_{\R^2}\nabla f\cdot\nabla g \ dx , \ \ f,g \in C_0^{\infty}(\R^2).
\]
Then as shown in this section, $(\E^n,C_0^{\infty}(\R^2))$ is closable in $L^2(\R^2,\rho dx)$ and its closure $(\E^n,D(\E^n))$ is the n-regularized Dirichlet form (see Theorem \ref{t;retrcos} (i)). Furthermore, $(\E^n,D(\E^n))$ is conservative (see Theorem \ref{t;retrcos} (i)).
Therefore, there exists a Hunt process $((\mathcal{B}^n_t)_{t \ge 0}, (\P_x)_{x \in \R^2_{\Delta}})$ associated with the Dirichlet form $(\E^n,D(\E^n))$  as follows (see Theorem \ref{th;ehplb} and \ref{t;lbmdf}):

\begin{thm}
$P$-a.s. it holds for all $x \in \R^2$
\[
\mathcal{B}^n_t = x + \int_0^t e^{-\frac{\gamma}{2} X_n(\mathcal{B}^n_s) + \frac{\gamma^2}{4} E[X_n(\mathcal{B}^n_s)^2]   } dB_s, \quad t \ge 0,
\]
$\P_x$-a.s. where $B$ is a standard 2-dimensional Brownian motion on $\R^2$.
\end{thm}

\section{Elliptic regularity results and singular diffusions associated with nonsymmetric degenerate forms}\label{s;nonsydema}
In this section we consider a non-symmetric Dirichlet form $(\E,D(\E))$ in $L^2(\R^d,m)$ as defined in \cite[I. Definition 4.5]{MR}. Let  $(T_t)_{t > 0}$ (resp. $(\hat{T}_t)_{t > 0}$) and $(G_{\alpha})_{\alpha > 0}$ (resp. $(\hat{G}_{\alpha})_{\alpha > 0}$ ) be the strongly continuous contraction $L^2(\R^d, m)$-semigroup (resp. cosemigroup) and resolvent (resp. coresolvent) associated to $(\E,D(\E))$ and $(L,D(L))$ (resp. $(\hat{L},D(\hat{L}))$) be the corresponding generator (resp. cogenerator) (see \cite[Diagram 3, p. 39]{MR} and see \cite[I. 1]{MR} for the definitions).  Then $(T_t)_{t>0}$ (resp. $(\hat{T}_t)_{t>0}$) and $(G_{\alpha})_{\alpha > 0}$ (resp. $(\hat{G}_{\alpha})_{\alpha > 0}$) are sub-Markovian (cf. \cite[I. Section 4]{MR}). Here an operator $S$ is called sub-Markovian if $0 \le f \le 1$, $f \in L^2(\R^d,m)$ implies $0 \le Sf \le 1$.
 Therefore analogous to symmetric case, $(T_t)_{t >0}$ (resp. $(G_{\alpha})_{\alpha > 0}$) can be extended to contraction semigroups (resp. resolvents) on $L^1(\R^d,m)$ and so $(T_t)_{t >0}$ (resp. $(G_{\alpha})_{\alpha > 0}$) restricted to $L^1(\R^d,m) \cap L^{\infty}(\R^d,m)$ can be extended to strongly continuous contraction semigroups (resp. contraction resolvents) on all $L^r(\R^d,m)$, $r \in [1,\infty)$. We denote the corresponding operator families again by $(T_t)_{t > 0}$ and $(G_{\alpha})_{\alpha > 0}$ and let $(L_r, D(L_r))$ be the corresponding generator on $L^r(\R^d,m)$. Furthermore by \cite[I. Corollary 2.21]{MR}, it holds that $(T_t)_{t>0}$ is analytic on $L^2(\R^d,m)$ and then by Stein interpolation (cf. e.g. \cite[Lecture 10, Theorem 10.8]{AV}) $(T_t)_{t>0}$ is also analytic semigroup on $L^r(\R^d,m)$ for all $r \in (1, \infty)$.\\
Suppose furthermore $(\E, D(\E))$ is a strictly quasi-regular, strongly local, non-symmetric Dirichlet form in $L^2(\R^d,m)$. Then there exists a Hunt process $\tilde{\bM} = (\tilde{\Omega}, \tilde{\F}, (\tilde{\F})_{t \ge 0}, \tilde{\zeta}, (\tilde{X}_t)_{t \ge 0}, (\tilde{\P}_x)_{x \in \R^d _{  \Delta } })$ (strictly properly) associated with $(\E, D(\E))$ (see \cite[V.2.13]{MR}). Consider the strict capacity Cap$_{\E}$ of the non-symmetric Dirichlet form $(\E,D(\E))$ as defined in \cite[V.2.1]{MR}, i.e. 
\[
\text{Cap}_{\E} =  \text{cap}_{1,\hat{G}_1\varphi}
\]
for some fixed $\varphi \in L^1(\R^d,m) \cap \mathcal{B}_b(\R^d)$, $0 < \varphi \le 1$. Then following \cite[Section 2]{RoShTr},
we may hence assume that 
\begin{equation}\label{contipath}
\tilde{\Omega} = \{\omega = (\omega (t))_{t \ge 0} \in C([0,\infty),\R^d_{\Delta}) \ | \ \omega(t) = \Delta \quad \forall t \ge \zeta(\omega) \}
\end{equation}
and
\[
\tilde{X}_t(\omega) = \omega(t), \quad t \ge 0.
\]
The existence of a Hunt process $\tilde{\bM}$ (associated with a strictly quasi-regular Dirichlet form $(\E,D(\E))$)  with \eqref{contipath} is used to construct another Hunt process we need in this section (see \cite[Theorem 2.12]{RoShTr}).

\subsection{Singular diffusions associated with non-symmetric degenerate forms}\label{s;ndee}
Throughout this subsection we assume:
\begin{itemize}
\item[(H1)]
$\rho = \xi^2$, $\xi \in H^{1,2}_{loc}(\R^d, dx)$,  $\rho > 0 \ \ dx$-a.e. and 
\[
\frac{\| \nabla \rho \|}{\rho} \in L^{p}_{loc} (\R^d, m), \quad m := \rho dx,
\]
$p:=(d + \varepsilon) \vee 2$ for some $\varepsilon >0$. 
\item[(H2)]  
Let $A=(a_{ij})_{1 \le i,j \le d}$ be a symmetric (possibly) degenerate (uniformly weighted)  elliptic $d \times d$ matrix, i.e. there exists a constant $\lambda \ge 1$
 such that for $dx$-a.e. $x \in \R^d$
\begin{equation}\label{eq;nosymde} 
\lambda^{-1} \ \rho(x) \ \| \xi \|^2 \le \langle A(x) \xi, \xi \rangle \le \lambda \ \rho(x) \ \|\xi\|^2, \quad  \forall \xi \in \R^d.
\end{equation} 
\end{itemize}
By (H1) the symmetric positive definite bilinear form
\begin{equation}\label{dibif}
\E^0(f,g) : = \frac{1}{2} \int_{\R^d} \langle \nabla f ,  \nabla g \rangle \ dm , \quad f, g \in C_0^{\infty}(\R^d)
\end{equation}
is closable in $L^2(\R^d,m)$ and its closure $(\E^0,D(\E^0))$ is a symmetric, strongly local, regular Dirichlet form (see \cite{AKR}). According to (H2) and \eqref{dibif}, the symmetric positive definite bilinear form
\[
\E^A(f,g) : = \frac{1}{2} \int_{\R^d} \langle A \nabla f ,  \nabla g \rangle \ dx, \quad f, g \in C_0^{\infty}(\R^d)
\]
is also closable in $L^2(\R^d,m)$ and its closure $(\E^A,D(\E^A))$ is a symmetric, strongly local, regular Dirichlet form (see \cite{FOT}). 
We further assume
\begin{itemize}
\item[(H3)] $B : \R^d \to \R^d, \ \|B\| \in L_{loc}^{p}(\R^d, m)$ where $p$ is the same as in (H1) and
\[
\int_{\R^d} \langle B,\nabla f \rangle \ dm = 0, \quad \forall f \in C_0^{\infty}(\R^d),
\]
\end{itemize}
and
\begin{itemize}
\item[(H4)]
\[
\left| \int_{\R^d} \langle B,\nabla f \rangle \ g \ \rho \ dx \right| \le c_0 \ \E^A_1 (f,f)^{1/2} \ \E^A_1 (g,g)^{1/2}, \quad \forall f,g \in C_0^{\infty}(\R^d),
\]
where $c_0$ is some constant (independent of $f$ and $g$) and $\E_{\alpha}^A(\cdot,\cdot):=\E^A(\cdot,\cdot) + \alpha (\cdot,\cdot)_{L^2(\R^d,m)}$, $\alpha > 0$.
\end{itemize}

Now we consider the non-symmetric bilinear form 
\[
\E(f,g) : = \frac{1}{2} \int_{\R^d} \langle  A \nabla f , \nabla g \rangle \ dx - \int_{\R^d} \langle B,\nabla f \rangle \ g \  dm, \quad
 f, g \in C_0^{\infty}(\R^d).
\]
Using the assumptions (H1)-(H4) and the closability of $(\E^A,C_0^{\infty}(\R^d))$ one can show that $(\E, C_0^{\infty}(\R^d))$ is closable in $L^2(\R^d,m)$ and the closure $(\E,D(\E))$ is a non-symmetric Dirichlet form (see \cite[II. 2. d) p. 48, 49]{MR} for the proof of Dirichlet form). Furthermore by \cite[V. Proposition 2.12 (ii)]{MR} $(\E, D(\E))$ is strictly quasi-regular.
\begin{lemma}\label{l;rhoctn} 
Assume (H1). Then $\rho$ is in $H^{1,p}_{loc}(\R^d,dx)$ and $\rho$ has a continuous $dx$-version in $C_{loc}^{1-d/p}(\R^d)$.
\end{lemma}
\proof
By \cite[Corollary 2.2]{AKR} the assumption (H1) implies $\rho \in H^{1,p}_{loc}(\R^d,dx)$. Therefore $\rho$ has a continuous $dx$-version in $C_{loc}^{1-d/p}(\R^d)$.
\qed

From now on, we shall always consider the continuous $dx$-version of $\rho$ as in Lemma \ref{l;rhoctn} and denote it also by $\rho$. We further assume the integrability of the derivative of $a_{ij}$:
\begin{itemize}
\item[(H5)]  For each $1 \le i,j \le d$, $a_{ij} \in H^{1,1}_{loc} (\R^d,dx)$ and $\partial_{j} \left( \frac{a_{ij}}{\rho} \right) \in L^p_{loc} (\R^d,dx)$, where $p$ is the same as in (H1). 
\end{itemize}
\begin{lemma}\label{l;integaio}
For each $1 \le i,j \le d$, $\frac{a_{ij}}{\rho} \in L^{\infty}(\R^d,m)$ and $\frac{\partial_{j} a_{ij}}{\rho} \in L^{p}_{loc}(\R^d,m)$.
\end{lemma}
\proof
By (H2), $ \left| \frac{a_{ij}}{\rho} \right |  \le  \lambda - \frac{1}{\lambda}(1-\delta_{ij})$, where $1 \le i,j \le d$ and $\delta_{ij}: = 1$ if $i=j$, $\delta_{ij}: = 0$ if $i \neq j$. Together with (H1) this implies that $\frac{a_{ij}}{\rho} \in L^{\infty}(\R^d,m)$. Furthermore by (H1), (H2), (H5)
\[
\frac{\partial_{j} a_{ij}}{\rho} =  \partial_{j} \left( \frac{a_{ij}}{\rho} \right) - \frac{a_{ij}}{\rho} \frac{\partial_j \rho}{\rho} \in L^{p}_{loc}(\R^d,m).
\]
\qed

By (H3) and Lemma \ref{l;integaio} we get                        
$ C_0^{\infty}(\R^d) \subset D(L_r)$ for any $r \in [1,p]$ and
\begin{equation}\label{eq;nsdge}
L _r \varphi =  \frac{1}{2}\sum_{i,j = 1}^{d}  \frac{a_{ij}}{\rho}  \  \partial_{ij} \varphi  +   \sum_{i = 1}^{d} \Big( \sum_{j = 1}^{d}   \frac{\partial_{j} a_{ij}}{2\rho}    + b_i \Big)   \ \partial_i \varphi, \quad \varphi \in C_0^{\infty}(\R^d), \quad r \in [1,p].
\end{equation}

The regularity properties of $(T_t)_{t > 0}$ and $(G_{\alpha})_{\alpha > 0}$ in the symmetric case (resp. non-symmetric case) were discussed in \cite{AKR}, \cite{BGS}, and Section \ref{symer} (resp. \cite{RoShTr}). Using Proposition \ref{p;gdegep} and 
Proposition \ref{t;morrey2.4} we obtain similar regularity properties for non-symmetric case with the (possibly) degenerate matrix $A$ satisfying (H2).
\begin{thm}\label{c;ressob}
Let $\alpha > 0$.
Suppose $g \in L^r(\R^d, m)$, $r \in [p,\infty)$. Then
\[
\rho  G_{\alpha} g \in H_{loc}^{1,p}(\R^d,dx)
\]
and for any open ball $B'\subset \overline{B'}\subset B \subset \overline{B}  \subset \{\rho > 0\}$ there exists $c_{B,\alpha} \in (0,\infty)$, independent of $g$, such that 
\begin{equation}\label{conresol}
\| \ \rho  G_{\alpha} g \ \|_{H^{1,p}(B',dx)} \le c_{B,\alpha} \ \Big(\|G_{\alpha} g \|_{L^1(B,m)} + \| g \|_{L^{p}(B,m)} \Big).
\end{equation}
\end{thm}
\proof
Let $g \in C_0^{\infty}(\R^d)$. Then we have by (H3) and  \eqref{eq;nsdge}
\[
\hat{L} \varphi =  \frac{1}{2}\sum_{i,j = 1}^{d}  \frac{a_{ij}}{\rho}  \  \partial_{ij} \varphi  +   \sum_{i = 1}^{d}  \Big( \sum_{j = 1}^{d}   \frac{\partial_{j} a_{ij}}{2 \rho}    - b_i \Big) \   \partial_i \varphi
\]
and
\begin{equation}\label{eq;geraseq}
\int_{\R^d} (\alpha - \hat{L})\varphi \ G_{\alpha} g \ \rho \ dx = \int_{\R^d} \varphi \ g \ \rho \ dx, \quad \forall \varphi \in C_0^{\infty}(\R^d).
\end{equation}
Now we apply Proposition \ref{p;gdegep} with $d_{ij} = \frac{a_{ij}}{2 \rho}$, $h_i =  \sum_{j = 1}^{d}   \frac{\partial_{j} a_{ij}}{2 \rho}    - b_i$, $c= -\alpha$, $\mu = - \rho G_{\alpha} g dx$, and $f=g \rho$. We check all necessary conditions: Since  $G_{\alpha} g \in L^1(\R^d,m) \cap L^{\infty} (\R^d,m)$ for $g \in C_0^{\infty}(\R^d) \subset L^1(\R^d,m) \cap L^{\infty}(\R^d,m)$, we can define the locally finite signed Borel measure $\mu = - \rho G_{\alpha}g dx$. By Lemma \ref{l;integaio} and (H5),  $ \frac{a_{ij}}{2 \rho} \in H^{1,p}_{loc}(\R^d,dx)$. By (H3) and Lemma \ref{l;integaio}, $\sum_{j = 1}^{d}   \frac{\partial_{j} a_{ij}}{2 \rho}    - b_i \in L^p_{loc}(\R^d,m)$ and so together with $G_{\alpha} g \in L^{\infty} (\R^d,m)$  this implies that $\sum_{j = 1}^{d}   \frac{\partial_{j} a_{ij}}{2 \rho}    - b_i \in L^p_{loc}(\R^d,\mu)$. Clearly $-\alpha \in L^p_{loc}(\R^d,\mu)$ and $g \rho \in L^p_{loc}(\R^d,dx)$. 

Hence for $g \in C_0^{\infty}(\R^d)$, $\rho G_{\alpha} g \in H_{loc}^{1,p}(\R^d,dx)$.  Let $u: =\rho G_{\alpha} g$. Using integration by parts \eqref{eq;geraseq} can be written as
\[
\int_{\R^d}  \sum_{i= 1}^{d}  \partial_{i} \varphi \left[  \sum_{j = 1}^{d}    \left( \frac{a_{ij}}{2 \rho} \right) \partial_j u    +  \left( \sum_{j = 1}^{d} \partial_j  \left( \frac{a_{ij}}{2 \rho} \right)     -  \sum_{j = 1}^{d}      \frac{\partial_{j} a_{ij}}{ 2 \rho}   +  b_i  \right)  u   \right] +  \varphi \ (\alpha u  -  \ \rho g )\ dx = 0.
\]
Therefore for any open balls $B'$, $B$ with $B'\subset \overline{B'}\subset B \subset \overline{B}  \subset \{\rho > 0\}$,  we can apply Proposition \ref{t;morrey2.4} with $\varphi \in C_0^{\infty}(B)$. We can check that all necessary conditions are satisfied. Therefore \eqref{conresol} holds for $g \in C_0^{\infty}(\R^d)$. Since $C_0^{\infty}(\R^d)$ is dense in ($L^r(\R^d, m)$, $\| \cdot \|_{L^r(\R^d, m)})$, $r \in [p,\infty) $, the assertion for general $g \in L^r(\R^d, m)$ follows by continuity.
\qed

\begin{cor}\label{c;dnsesob}
Let $t>0$, $r \in [p,\infty)$. 
\begin{itemize}
\item[(i)] Let $u \in D(L_r)$. Then 
\[
\rho T_t u \in H^{1,p}_{loc} (\R^d, dx)
\]
and for any open balls $B'$, $B$ with  $B'\subset \overline{B'}\subset B \subset \overline{B} \subset  \{\rho > 0\}$ there exists $c_B \in (0,\infty)$ (independent of $u$ and $t$) such that
\begin{equation}\label{eq;nosmdej}
\| T_t u\|_{H^{1,p}(B',dx)} \le c_B \left( m(B)^{\frac{r-1}{r}} \| u \|_{L^r (\R^d,m) }  + m(B)^{\frac{r-p}{rp}} \| (1-L_r) u \|_{L^r(\R^d,m)}    \right).
\end{equation}
\item[(ii)] Let $u \in L^r(\R^d,m)$. Then the above statements still hold with \eqref{eq;nosmdej} replaced by
\[
\|\rho T_t u\|_{H^{1,p}(B',dx)} \le  \tilde{c}_B \left(1+ \frac{1}{t}\right)  \| u \|_{L^r(\R^d,m)},
\]
where $\tilde{c}_B \in (0,\infty)$ (independent of $f$, $t$).
\end{itemize}

\end{cor}
\proof
Using Theorem \ref{c;ressob} and the analyticity of $(T_t)_{t>0}$ the proof is the same as \cite[Corollary 2.4]{AKR}.
\qed

In this subsection we shall keep the notation
\[
E : = \{\rho > 0 \}.
\] 
By virtue of Lemma \ref{l;rhoctn}, Corollaries \ref{c;ressob}, \ref{c;dnsesob}, exactly as in \cite[Section 3]{AKR} (cf. \cite{RoShTr}),
we obtain the existence of a transition kernel density $p_t(\cdot,\cdot)$ on the open set $E$
such that
\[
P_t f(x) : = \int_{\R^d} f(y) p_t(x,y) \ m(dy), \quad x \in E, \ t > 0
\]
is a submarkovian transition function and an $m$-version of $T_t f $ for any $f \in \cup_{r \ge p} L^r(\R^d,m)$. Moreover, $(P_t)_{t >0}$ satisfies the properties as in Theorem \ref{t;2.7gp} with $\R^d$ replaced by $E$.\\
Let Cap (resp. Cap$_0$) be the capacity related to the symmetric Dirichlet form ($\E^A,D(\E^A)$) (resp. ($\E^0,D(\E^0)$)) as defined in \cite[Section 2.1]{FOT}. Note that Cap$_0(\{\rho=0\})=0$ by \cite[Theorem 2]{fuku85}. 
\begin{lemma}\label{l;ca3im}
Let $N \subset \R^d$. Then 
\[
\emph{Cap}_0(N) = 0  \Rightarrow \emph{Cap}(N) = 0  \Rightarrow \emph{Cap}_{\E}(N) = 0.
\]
In particular $\emph{Cap}_{\E}(\{\rho = 0\})=0$.
\end{lemma}
\proof
Let $N \subset \R^d$ be such that Cap$_0(N)=0$. By \eqref{eq;nosymde}, Cap$(N)=0$. Then the statement follows from the proof of 
\cite[Lemma 2.10]{RoShTr} that Cap$_{\E}(N) = 0$.
\qed

Following \cite[Section 2]{RoShTr} with the existence of $\tilde{\bM}$, Theorem \ref{c;ressob}, Corollary \ref{c;dnsesob}, and Lemma \ref{l;ca3im}, we obtain (cf. \cite[Theorem 2.12]{RoShTr}):
\begin{thm}\label{th;ndexisthunt}
There exists a continuous Hunt process
\[
\bM =  (\Omega, \F, (\F_t)_{t \ge 0}, \zeta, (X_t)_{t \ge 0}, (\P_x)_{x \in E_{\Delta}}   )
\]
with state space $E$, having the transition function $(P_t)_{t \ge 0}$ as transition semigroup.
\end{thm}

We further consider 
\begin{itemize}
\item[(H6)] $(\E,D(\E))$ is conservative.
\end{itemize}

\begin{remark}\label{r;cnsenode}
\begin{itemize}
\item[(i)] Assume (H6) holds (additionally to (H1)-(H5)). Then Theorem \ref{t;trsyst} and Remark \ref{r;sytr} hold with $\R^d$ replaced by $E$.
Furthermore, one can drop $\Delta$ in Theorem \ref{th;ndexisthunt}
\item[(ii)] Assume (H1)-(H5). Then Lemma \ref{l;trsydepsr} and Lemma \ref{potwitrstr} hold with $\R^d$ replaced by $E$.
\end{itemize}
\end{remark}

\begin{lemma}\label{l;gesm2u}
\begin{itemize}
\item[(i)]
For $u \in C_0^{\infty} (\R^d)$
\[
L u^2 - 2 u \ Lu = \sum_{i,j = 1}^{d}  \frac{a_{ij}}{\rho}  \ \partial_{i} u \  \partial_{j} u.
\]
\item[(ii)]
Let $u \in C_0^{\infty}(\R^d)$ and
\[
M_t : = \left( u(X_t) - u(X_0) - \int_0^t Lu(X_r) \ dr \right)^2 - \int_0^t  \Big(   \sum_{i,j = 1}^{d}  \frac{a_{ij}}{\rho}  \ \partial_{i} u \  \partial_{j} u \Big) (X_r) \ dr, \quad t \ge 0.
\]
Then $(M_t)_{t \ge 0}$ is  an $(\mathcal{F}_t)_{t \ge 0}$-martingale under $\P_x$, $\forall x \in E$.
\end{itemize}
\end{lemma}
\proof
(i) This follows immediately from \eqref{eq;nsdge}.  (ii) Using (i) and Remark \ref{r;cnsenode} (ii) the proof is similar to \cite[Proposition 3.3]{RoShTr}.
\qed

Let 
\begin{equation}\label{e;mademp}
M_t^{u} : = u(X_t) - u(X_0) - \int_0^t Lu(X_s) \ ds, \quad u \in C_0^{\infty}(E), \quad \ t \ge 0 .
\end{equation}
Clearly $(M^{u}_{t})_{t \ge 0}$ is a continuous $(\mathcal{F}_t)_{t \ge 0}$-martingale under $\P_x$, $x \in E$ and for $u_1$, $u_2  \in C_0^{\infty}(E)$, $M_t^{u_1 + u_2} = M_t^{u_1} + M_t^{u_2}$. By Lemma \ref{l;gesm2u} $M_t^{u} \in L^2(\Omega, \F_t, \P_x)$ and  its quadratic variation is given by
\[
\langle M^{u} \rangle_t =  \int_0^t  \Big(   \sum_{i,j = 1}^{d}  \frac{a_{ij}}{\rho}  \ \partial_{i} u \  \partial_{j} u \Big) (X_s) \ ds.
\]  
We define quadratic covariation process by
\[
\langle M^{u_1}, M^{u_2} \rangle = \frac{1}{2} \Big(  \langle M^{u_1} + M^{u_2} \rangle  - \langle M^{u_1} \rangle  -\langle M^{u_2} \rangle    \Big).
\]
\begin{lemma}\label{le;quadco}
Let  $u_1$, $u_2 \in C_0^{\infty}(E)$.
Then
\[
\langle M^{u_1}, M^{u_2} \rangle_t = \int_0^t \Big(   \sum_{i,j = 1}^{d}  \frac{a_{ij}}{\rho}  \ \partial_{i} u_1 \  \partial_{j} u_2 \Big) (X_s) \ ds.
\]
\end{lemma}

\begin{thm}\label{t;fdnons}
Under (H1)-(H5) after enlarging the stochastic basis $(\Omega, \F, (\F_t)_{t\ge 0},\P_x )$ appropriately for any $x\in E$, $i=1,\dots,d$, the process $\bM$ satisfies 
\[
X_t^i = x_i + \sum_{j=1}^d \int_0^t \frac{\sigma_{ij}}{\sqrt{\rho}} (X_s) \ dW_s^j +  \int^{t}_{0}   \left(\sum_{j=1}^d  \frac{ \partial_j a_{ij}}{2 \rho} + b_i    \right) (X_s) \, ds, \quad t < \zeta,
\]
where $(\sigma_{ij})_{1 \le i,j \le d} =  \sqrt{A} $ is the positive square root of the matrix $A$, $W = (W^1,\dots,W^d)$ is a standard d-dimensional Brownian motion on $\R^d$. If additionally (H6) holds, then $\zeta$ can be replaced by $\infty$ (cf. Remark \ref{r;cnsenode} (i)).
\end{thm}
\proof
Let $u_k \in C_0^{\infty}(E)$, $k = 1,\dots,d$. Note that by Lemma \ref{le;quadco}
\[
\langle M^{u_k}, M^{u_l} \rangle_t = \int_0^t \Big(   \sum_{i,j = 1}^{d}  \frac{a_{ij}}{\rho}  \ \partial_{i} u_k \  \partial_{j} u_l \Big) (X_s) \ ds, \quad 1 \le k,l  \le d.
\]

Suppose $\zeta<\infty$. Since $\Big(   \sum_{i,j = 1}^{d}  \frac{a_{ij}}{\rho}  \ \partial_{i} u_k \  \partial_{j} u_l \Big)$  is degenerate, it is standard that  there is an enlargement $(\bar{\Omega}, \bar{\F}, \bar{\P}_x )$ of the underlying probability space $(\Omega, \F, \P_x )$ 
and a d-dimensional Brownian motion $(W_t)_{t \ge 0} = (W_t^1,\dots, W_t^d)_{t \ge 0}$ on $(\bar{\Omega}, \bar{\F}, \bar{\P}_x )$ and a $d \times d$ matrix $\eta = (\eta_{ki})_{1 \le i,k \le d}$ such that
\[
M_t^{u_k} = \sum_{i=1}^{d}  \int_0^t \ \eta_{ki} (X_s) \ dW_s^k, \quad 1 \le k \le d, \  t\ge 0.
\]
Here  $\langle  \rho^{-1} A  \nabla u_k, \nabla u_l \rangle = \langle  \sqrt{\rho}^{-1} \sqrt{A}  \nabla u_k, \sqrt{\rho}^{-1} \sqrt{A} \nabla u_l \rangle =  \sum_{i=1}^{d} \eta_{ki} \eta_{li}$ (cf. \cite[Section 3.4.A., 4.2 Theorem]{KS}).
The identification of $X$ up to $\zeta$ is now obtained by using Remark \ref{r;cnsenode} (ii) with an appropriate localizing sequence for which the coordinate projections on $E$ coincide locally with $C_0^{\infty}(E)$-functions.  The localization of the drift part in \eqref{e;mademp} is trivial.
\qed

\begin{remark}\label{r;ensdb}
The results of this subsection include the particular case with
\[
A(x) = \rho (x) \cdot Id \quad \text{and} \ \lambda =1.  
\]
Therefore, our consequences in this subsection lead to an extension of the results of \cite{RoShTr}.
\end{remark}

\subsection{Singular diffusions associated with locally uniformly elliptic form}\label{s;3.2sdlce}
In this subsection, we alway assume (H1) for $\rho$ and (H3) for $B$ as in Subsection \ref{s;ndee}  and consider a symmetric locally uniformly elliptic $d \times d$ matrix $A=(a_{ij})_{1 \le i,j \le d}$ such that 
\begin{itemize}
\item[(H7)]  
$a_{ij}\in H^{1,p}_{loc} (\R^d,dx)$, $p$ as in (H1) and for any compact set $K \subset \R^d$ there exists a constant $c_K >0$
 such that for $dx$-a.e. $x \in K$
\[
c_K  \ \| \xi \|^2 \le \langle A(x) \xi, \xi \rangle, \quad  \forall \xi \in \R^d.
\]
\end{itemize}
By Lemma \ref{l;rhoctn} (resp. (H7)), $\rho$ (resp. $a_{ij}$) has a continuous $dx$-version in $C_{loc}^{1- d/p} (\R^d)$.
Throughout this subsection we shall always consider the continuous $dx$-version of $\rho$ (resp. $a_{ij}$)  and denote it also by $\rho$ (resp. $a_{ij}$). The symmetric positive definite bilinear form
\begin{equation}\label{e;dafo3.11}
\E^A(f,g) : = \frac{1}{2} \int_{\R^d} \langle A \nabla f ,  \nabla g \rangle \ dm, \quad f, g \in C_0^{\infty}(\R^d), \quad m: = \rho dx
\end{equation}
is  closable in $L^2(\R^d,m)$ and its closure $(\E^A,D(\E^A))$ is a symmetric, strongly local, regular Dirichlet form (see \cite[II. Exercise 2.4]{MR}). 
We further assume (H4) with $\E^A$ replaced by $\E^A$ as in \eqref{e;dafo3.11}. Now we consider the non-symmetric bilinear form 
\[
\E(f,g) : = \frac{1}{2} \int_{\R^d} \langle  A \nabla f , \nabla g \rangle \ dm - \int_{\R^d} \langle B,\nabla f \rangle \ g \  dm, \quad
 f, g \in C_0^{\infty}(\R^d).
\]
Since $(\E^A,C_0^{\infty}(\R^d))$ is closable, $(\E, C_0^{\infty}(\R^d))$ is also closable in $L^2(\R^d,m)$ and the closure $(\E,D(\E))$ is a non-symmetric Dirichlet form (see \cite[II. 2. d) p. 48, 49]{MR}). Furthermore by \cite[V. Proposition 2.12 (ii)]{MR} $(\E, D(\E))$ is strictly quasi-regular. By (H1), (H3), and (H7), we get $C_0^{\infty}(\R^d) \subset  D(L_r)$ and
\begin{equation}\label{noloung}
L_r f =  \frac{1}{2}\sum_{i,j = 1}^{d}  a_{ij}  \  \partial_{ij} f  +   \sum_{i = 1}^{d}  \Big( \sum_{j = 1}^{d}   \frac{1}{2} \partial_{j} a_{ij} + \frac{\partial_j \rho }{2 \rho} a_{ij}   + b_i \Big) \   \partial_i f, \quad f \in C_0^{\infty}(\R^d), \quad r \in [1,p],
\end{equation}
where $p$ is same as in (H1).

\begin{thm}
Let $\alpha > 0$.
Suppose $g \in L^r(\R^d, m)$, $r \in [p,\infty)$. Then
\[
\rho  G_{\alpha} g \in H_{loc}^{1,p}(\R^d,dx)
\]
and for any open ball $B'\subset \overline{B'}\subset B \subset \overline{B}  \subset \{\rho > 0\}$ there exists $c_{B,\alpha} \in (0,\infty)$, independent of $g$, such that 
\begin{equation}\label{eq;esglue}
\| \ \rho  G_{\alpha} g \ \|_{H^{1,p}(B',dx)} \le c_{B,\alpha} \ \Big(\|G_{\alpha} g \|_{L^1(B,m)} + \| g \|_{L^{p}(B,m)} \Big).
\end{equation}
\end{thm}
\proof
For $g \in C_0^{\infty}(\R^d)$ we have
\begin{equation}\label{eq;geraseque}
\int (\alpha - \hat{L})\varphi \ G_{\alpha} g \ \rho \ dx = \int \varphi \ g \ \rho \ dx, \quad \forall \varphi \in C_0^{\infty}(\R^d),
\end{equation}
where
\[
\hat{L} \varphi =  \frac{1}{2}\sum_{i,j = 1}^{d}  a_{ij}  \  \partial_{ij} \varphi  +   \sum_{i = 1}^{d}  \Big( \sum_{j = 1}^{d}   \frac{1}{2} \partial_{j} a_{ij} + \frac{\partial_j \rho }{2 \rho} a_{ij}   - b_i \Big) \   \partial_i \varphi.
\]
Now we apply Proposition \ref{p;gdegep} with $d_{ij} = \frac{a_{ij}}{2 }$, $h_i =  \sum_{j = 1}^{d}   \frac{1}{2} \partial_{j} a_{ij} + \frac{\partial_j \rho }{2 \rho} a_{ij}   - b_i$, $c= -\alpha$, $\mu = - \rho G_{\alpha} g dx$, and $f=g \rho$ to prove the assertion for $g \in C_0^{\infty}(\R^d)$. All necessary integrability conditions are satisfied. Hence for $g \in C_0^{\infty}(\R^d)$, $\rho G_{\alpha} g \in H_{loc}^{1,p}(\R^d,dx)$.  Let $u: =\rho G_{\alpha} g$. Using integration by parts \eqref{eq;geraseque} can be written as
\[
\int_{\R^d}  \sum_{i= 1}^{d}  \partial_{i} \varphi \left[  \sum_{j = 1}^{d}    \left( \frac{a_{ij}}{2 } \right) \partial_j u    +  \left( \sum_{j = 1}^{d} \partial_j  \left( \frac{a_{ij}}{2 } \right)     -  \sum_{j = 1}^{d}   \frac{1}{2} \partial_{j} a_{ij} - \frac{\partial_j \rho }{2 \rho} a_{ij}   + b_i   \right)  u   \right] +  \varphi \ (\alpha u  -  \ \rho g )\ dx = 0.
\]
Therefore for any open balls  $B'$, $B$ with $B'\subset \overline{B'}\subset B \subset \overline{B}  \subset \{\rho > 0\}$, we can apply Proposition \ref{t;morrey2.4} with $\varphi \in C_0^{\infty}(B)$. Since all integrability conditions are satisfied, \eqref{eq;esglue} holds for $g \in C_0^{\infty}(\R^d)$. Then, since $C_0^{\infty}(\R^d)$ is dense in ($L^r(\R^d, m)$, $\| \cdot \|_{L^r(\R^d, m)})$, $r \in [p,\infty) $, the assertion for general $g \in L^r(\R^d, m)$ follows by continuity.
\qed

Let  $E:=\{ \rho > 0 \}$. Then we obtain Corollary \ref{c;dnsesob} and following subsequent results as in Section \ref{s;ndee}, 
we have the existence of a transition kernel density $p_t(\cdot,\cdot)$ on the open set $E$
such that
\[
P_t f(x) : = \int_{\R^d} f(y) p_t(x,y) \ m(dy), \quad x \in E, \ t > 0
\]
is a submarkovian transition function and an $m$-version of $T_t f $ for any $f \in \cup_{r \ge p} L^r(\R^d,m)$.\\
Let Cap be the capacity related to the symmetric Dirichlet form ($\E^A,D(\E^A)$) as defined in \cite[Section 2.1]{FOT}.
\begin{lemma}
\begin{itemize}
\item[(i)] $\emph{Cap}(\{\rho = 0\})=0$.
\item[(ii)]
Let $N \subset \R^d$. Then 
\[
\emph{Cap}(N) = 0  \Rightarrow \emph{Cap}_{\E}(N) = 0.
\]
In particular $\emph{Cap}_{\E}(\{\rho = 0\})=0$.
\end{itemize}
\end{lemma}
\proof
(i) We follow the proof of \cite[Proposition 4.4]{BGS} (cf. \cite{fuku85}). (ii) The proof is the same as Lemma \ref{l;ca3im}.
\qed

Then analogous to Theorem \ref{th;ndexisthunt}, we can construct a continuous Hunt process
\[
\bM =  (\Omega, \F, (\F_t)_{t \ge 0}, \zeta, (X_t)_{t \ge 0}, (\P_x)_{x \in E_{\Delta}}   )
\]
with state space $E$, having the transition function $(P_t)_{t \ge 0}$ as transition semigroup (for details see Subsection \ref{s;ndee}).
\begin{lemma}\label{l;3.19gr}
\begin{itemize}
\item[(i)] For $u \in C_0^{\infty} (\R^d)$
\[
L u^2 - 2 u \ Lu = \sum_{i,j = 1}^{d}  a_{ij} \ \partial_{i} u \  \partial_{j} u.
\]
\item[(ii)] Let $u \in C_0^{\infty}(\R^d)$ and
\[
M_t : = \left( u(X_t) - u(X_0) - \int_0^t Lu(X_r) \ dr \right)^2 - \int_0^t  \Big(   \sum_{i,j = 1}^{d}  a_{ij}  \ \partial_{i} u \  \partial_{j} u \Big) (X_r) \ dr, \quad t \ge 0.
\]
Then $(M_t)_{t \ge 0}$ is  an $(\mathcal{F}_t)_{t \ge 0}$-martingale under $\P_x$, $\forall x \in E$.
\end{itemize}
\end{lemma}
\proof
(i) This follows immediately from \eqref{noloung}. (ii) For the proof, we refer to Lemma \ref{l;gesm2u}.
\qed

\begin{thm}
Assume (H1), (H3), (H4), and (H7). For each $x\in E$, $i=1,\dots,d$, the process $\bM$ satisfies 
\[
X_t^i = x_i + \sum_{j=1}^d \int_0^t  \sigma_{ij} (X_s) \ dW_s^j +  \int^{t}_{0}   \left(\sum_{j=1}^d  \frac{ \partial_j a_{ij}}{2} +  \frac{\partial_j \rho }{2 \rho} a_{ij} + b_i    \right) (X_s) \, ds, \quad t < \zeta,
\]
where $(\sigma_{ij})_{1 \le i,j \le d} =  \sqrt{A} $ is the positive square root of the matrix $A$, $W = (W^1,\dots,W^d)$ is a standard d-dimensional Brownian motion on $\R^d$. If we additionally assume conservativeness of $(\E,D(\E))$, then $\zeta$ can be replaced by $\infty$.
\end{thm}
\proof
Using Lemma \ref{l;3.19gr} and previous results, the proof is similar to Theorem \ref{t;fdnons}. So we omit it.
\qed

\subsection{Singular diffusions associated with degenerate forms with Lebesgue measure}\label{s;3.3lms}
In this subsection we consider the following assumption:
\begin{itemize}
\item[(H8)]
Let $\psi \in C(\R^d)$ and $\psi > 0$ $dx$-a.e. and $A=(a_{ij})_{1 \le i,j \le d}$ be a symmetric degenerate  elliptic $d \times d$ matrix, that is  $a_{ij} \in H^{1,p}_{loc}(\R^d,dx)$, $p > d$, such that for $dx$-a.e. $x \in \R^d$
\[
\psi (x) \ \| \xi \|^2  \le \langle A(x) \xi, \xi \rangle, \quad  \forall \xi \in \R^d.
\]
\end{itemize}

By \cite[Section 3.1 (1$^{\circ}$)]{FOT} the symmetric positive definite bilinear form
\[
\E^A(f,g) : = \frac{1}{2} \int_{\R^d} \langle A \nabla f ,  \nabla g \rangle \ dx, \quad f, g \in C_0^{\infty}(\R^d)
\]
is  closable in $L^2(\R^d,dx)$ and its closure $(\E^A,D(\E^A))$ is a symmetric, strongly local, regular Dirichlet form. We further assume:

\begin{itemize}
\item[(H9)] $B : \R^d \to \R^d, \ \|B\| \in L_{loc}^{p}(\R^d, dx)$ where $p > d$ and
\[
\int_{\R^d} \langle B,\nabla f \rangle \ dx = 0, \quad \forall f \in C_0^{\infty}(\R^d),
\]
and
\[
\left| \int_{\R^d} \langle B,\nabla f \rangle \ g \ dx \right| \le c_0 \ \E^A_1 (f,f)^{1/2} \ \E^A_1 (g,g)^{1/2}, \quad \forall f,g \in C_0^{\infty}(\R^d),
\]
where $c_0$ is some constant (independent of $f$ and $g$).
\end{itemize}

Now we consider the non-symmetric bilinear form 
\[
\E(f,g) : = \frac{1}{2} \int_{\R^d} \langle  A \nabla f , \nabla g \rangle \ dx - \int_{\R^d} \langle B,\nabla f \rangle \ g \  dx, \quad
 f, g \in C_0^{\infty}(\R^d).
\]
Since $(\E^A,C_0^{\infty}(\R^d))$ is closable, $(\E, C_0^{\infty}(\R^d))$ is also closable in $L^2(\R^d,dx)$ and its closure $(\E,D(\E))$ is a non-symmetric Dirichlet form (cf. \cite[II. 2. d)]{MR}). Furthermore by \cite[V. Proposition 2.12 (ii)]{MR} $(\E, D(\E))$ is strictly quasi-regular.
\begin{thm}
Let $\alpha > 0$.
Suppose $g \in L^r(\R^d, dx)$, $r \in [p,\infty)$. Then
\[
G_{\alpha} g \in H_{loc}^{1,p}(\R^d,dx)
\]
and for any open ball $B'\subset \overline{B'}\subset B \subset \overline{B}  \subset \{\psi > 0\}$ there exists $c_{B,\alpha} \in (0,\infty)$, independent of $g$, such that 
\begin{equation}\label{eq;relees}
\| G_{\alpha} g \ \|_{H^{1,p}(B',dx)} \le c_{B,\alpha} \ \Big(\|G_{\alpha} g \|_{L^1(B,dx)} + \| g \|_{L^{p}(B,dx)} \Big).
\end{equation}
\end{thm}
\proof
Let $g \in C_0^{\infty}(\R^d)$. Then we have
\begin{equation}\label{eq;cogef}
\int (\alpha - \hat{L})\varphi \ G_{\alpha} g  \ dx = \int \varphi \ g   \ dx, \quad \forall \varphi \in C_0^{\infty}(\R^d),
\end{equation}
where
\[
\hat{L} \varphi =  \frac{1}{2}\sum_{i,j = 1}^{d}  a_{ij}  \  \partial_{ij} \varphi  +   \sum_{i = 1}^{d}  \Big( \sum_{j = 1}^{d}   \frac{\partial_{j} a_{ij}}{2}    - b_i \Big) \   \partial_i \varphi.
\]
Now we apply Proposition \ref{p;gdegep} with $d_{ij} = \frac{a_{ij}}{2}$, $h_i =  \sum_{j = 1}^{d}   \frac{\partial_{j} a_{ij}}{2 }    - b_i$, $c= -\alpha$, $\mu = - G_{\alpha} g dx$, and $f=g $ to prove the assertion for $g \in C_0^{\infty}(\R^d)$. All necessary integrability conditions are satisfied. Hence for $g \in C_0^{\infty}(\R^d)$, $G_{\alpha} g \in H_{loc}^{1,p}(\R^d,dx)$.  Let $u: = G_{\alpha} g$. Using integration by parts \eqref{eq;cogef} can be written as
\[
\int_{\R^d}  \sum_{i= 1}^{d}  \partial_{i} \varphi \left[  \sum_{j = 1}^{d}    \left( \frac{a_{ij}}{2} \right) \partial_j u    +  \left( \sum_{j = 1}^{d} \partial_j  \left( \frac{a_{ij}}{2} \right)     -  \sum_{j = 1}^{d}      \frac{\partial_{j} a_{ij}}{ 2 }   +  b_i  \right)  u   \right] +  \varphi \ (\alpha u  -  g )\ dx = 0.
\]
Now for any open balls $B'$, $B$ with $B'\subset \overline{B'}\subset B \subset \overline{B}  \subset \{\psi > 0\}$, we can apply Proposition \ref{t;morrey2.4} with $\varphi \in C_0^{\infty}(B)$.  Then, since all integrability conditions are satisfied, \eqref{eq;relees} holds for $g \in C_0^{\infty}(\R^d)$.
Since $C_0^{\infty}(\R^d)$ is dense in ($L^r(\R^d, dx)$, $\| \cdot \|_{L^r(\R^d, dx)})$, $r \in [p,\infty) $, the assertion for general $g \in L^r(\R^d, dx)$ follows by continuity.
\qed

Let $E_1 : = \{\psi > 0 \}$. Then we obtain Corollary \ref{c;dnsesob} with $\{\rho > 0\}$ replaced by $E_1$ and following subsequent results as in Section \ref{s;ndee}, 
we have the existence of a transition kernel density $p_t(\cdot,\cdot)$ on the open set $E_1$
such that
\[
P_t f(x) : = \int_{\R^d} f(y) p_t(x,y) \ dy, \quad x \in E_1, \ t > 0
\]
is a submarkovian transition function  and an $dx$-version of $T_t f $ for any $f \in \cup_{r \ge p} L^r(\R^d,dx)$. We further assume
\begin{itemize}
\item[(H10)] Cap$_{\E}(\{\psi = 0\}$)=0.
\end{itemize}
\begin{remark}
The assumption (H10) is satisfied if we take $\psi$ as in (H1) (see Subsection \ref{s;ndee}).
\end{remark}

Then, similar to Theorem \ref{th;ndexisthunt}, we can construct a continuous Hunt process
\[
\bM =  (\Omega, \F, (\F_t)_{t \ge 0}, \zeta, (X_t)_{t \ge 0}, (\P_x)_{x \in {E_{1}}_{\Delta}}   )
\]
with state space $E_1$, having the transition function $(P_t)_{t \ge 0}$ as transition semigroup (for details see Subsection \ref{s;ndee}).\\ 
By (H8) and (H9) we get                        
$ C_0^{\infty}(\R^d) \subset D(L_r)$ for any $r \in [1,p]$ and
\begin{equation}\label{eq;ch4genlm}
L _r f =  \frac{1}{2}\sum_{i,j = 1}^{d}  a_{ij}  \  \partial_{ij} f  +   \sum_{i = 1}^{d} \Big( \sum_{j = 1}^{d}   \frac{\partial_{j} a_{ij}}{2}    + b_i \Big)   \ \partial_i f , \quad f \in C_0^{\infty}(\R^d).
\end{equation}

\begin{lemma}\label{l;3.23len}
\begin{itemize}
\item[(i)]
For $u \in C_0^{\infty} (\R^d)$
\[
L u^2 - 2 u \ Lu = \sum_{i,j = 1}^{d}  a_{ij} \ \partial_{i} u \  \partial_{j} u.
\]
\item[(ii)]
Let $u \in C_0^{\infty}(\R^d)$ and
\[
M_t : = \left( u(X_t) - u(X_0) - \int_0^t Lu(X_r) \ dr \right)^2 - \int_0^t  \Big(   \sum_{i,j = 1}^{d}  a_{ij}  \ \partial_{i} u \  \partial_{j} u \Big) (X_r) \ dr, \quad t \ge 0.
\]
Then $(M_t)_{t \ge 0}$ is  an $(\mathcal{F}_t)_{t \ge 0}$-martingale under $\P_x$, $\forall x \in E_1$.
\end{itemize}
\end{lemma}
\proof
(i) This follows immediately from \eqref{eq;ch4genlm}. (ii) For the proof, we refer to Lemma \ref{l;gesm2u}.
\qed

\begin{thm}
Under (H8)-(H10) for any $x\in E_1$, $i=1,\dots,d$, the process $\bM$ satisfies 
\begin{equation*}
X_t^i = x_i + \sum_{j=1}^d \int_0^t  \sigma_{ij} (X_s) \ dW_s^j +  \int^{t}_{0}   \left(\sum_{j=1}^d  \frac{ \partial_j a_{ij}}{2} + b_i    \right) (X_s) \, ds, \quad t < \zeta,
\end{equation*}
where $(\sigma_{ij})_{1 \le i,j \le d} =  \sqrt{A} $ is the positive square root of the matrix $A$, $W = (W^1,\dots,W^d)$ is a standard d-dimensional Brownian motion on $\R^d$.  If we additionally assume conservativeness of $(\E,D(\E))$, then $\zeta$ can be replaced by $\infty$.
\end{thm}
\proof
Using Lemma \ref{l;3.23len} and previous results, the proof is similar to Theorem \ref{t;fdnons}.
\qed

\vspace*{2cm}
\noindent Jiyong Shin\\
School of Mathematics \\
Korea Institute for Advanced Study\\
85 Hoegiro Dongdaemun-gu,\\
Seoul 02445, South Korea,  \\
E-mail: yonshin2@kias.re.kr \\
\end{document}